\documentclass[oneside, 10pt]{article}
\usepackage{amsmath, graphicx, enumerate, xurl, tabularx, multirow}
\usepackage[utf8]{inputenc}
\usepackage{mathptmx}
\usepackage{amssymb} % For using \mathbb command
\usepackage[margin=4.5cm, textheight=18cm, top=6cm]{geometry}
\makeatletter
\DeclareFontFamily{U}{tipa}{}
\DeclareFontShape{U}{tipa}{m}{n}{<->tipa10}{}
\newcommand{\arc@char}{{\usefont{U}{tipa}{m}{n}\symbol{62}}}%
\newcommand{\arc}[1]{\mathpalette\arc@arc{#1}}

\newcommand{\arc@arc}[2]{%
  \sbox0{$\m@th#1#2$}%
  \vbox{
    \hbox{\resizebox{\wd0}{\height}{\arc@char}}
    \nointerlineskip
    \box0
  }%
}
\makeatother
\title{On the computation of the arcsin function in the Kerala school of astronomy and mathematics}
\author{{\bf Dr. V. N. Krishnachandran}\\
Principal (retd), Govt. Victoria College, Palakkad -- 678001\\
(Formerly) Professor of Mathematics,\\ 
Govt. Engineering College, Thrissur -- 680009\\
(email: {\tt krishnachandranvn@gmail.com})}
\date{}
\begin{document}
\maketitle
\tableofcontents
\newpage
\begin{abstract}
This paper  examines how the mathematicians and astronomers of the Kerala school tackled the  problem of computing the values of the arcsin function. Four different approaches are discussed all of which are found in N\=\i laka\d n\d tha Somay\=aj\=\i's (1444--1545 CE) {\em Tantrasa\.ngraha} and the roots of  all of which can be traced to ideas originally articulated by Sa\.ngamagr\=ama M\=adhava (c. 1340–1425 CE): (i) a simple method when the argument is small; (ii) an iterative method when the argument is small; (iii) a method based on a lookup table; (iv) a method when the argument is large. The paper also contains the original Sanskrit verses describing the various methods and English translations thereof.  Moreover, there is  a presentation of a novel method  for computing the circumference of a circle found in Jye\d s\d thadeva's (c. 1500--1575 CE) {\em Yuktibh\=a\d s\=a} which is based on method (i) for computing the arcsin function. All methods have been illustrated with numerical examples. A surprising by-product of the investigation is a totally unexpected appearance of a core integer sequence, namely, the entry A001764 in the Online Encyclopedia of Integer Sequence, while studying the iterative method for computing the arcsin function. 
\qquad\\[1mm]
\end{abstract}
\section{Introduction}
Sa\.ngamagr\=ama  M\=adhava's  (c. 1340–1425 CE) methods for computing the values of the jy\=a, kojy\=a and \'sara functions are well known. 
His methods to compute the arc length when its tangent is given is also well known. In this paper we examine how the mathematicians and astronomers of the Kerala school tackled the related problem of computing the values of the the arc lengths when the values of the jy\=a-s are given, or, equivalently, computing the values of the arcsin function.  
This problem is highly significant and it arises in connection with several astronomical problems especially in connection with problems involving the ``equation of centre'' of the various planets. 
In particular, the problem of computing values of the arcsin function arises in the computation of the {\em candrav\=akya}-s. Since M\=adhava is famous as the author of a more accurate set of {\em candrav\=akya}-s than the then available Vararuci's {\em candrav\=akya}-s, he definitely should have addressed the problem. 
No direct account of how M\=adhava had tackled the problem have survived. However, in N\=\i laka\d n\d tha Somay\=aj\=\i's (1444--1545 CE) {\em Tantrasa\.ngraha} there is a detailed discussion of the problem and it could be the case that the solution presented by Somay\=aj\=\i\ might have had its origins in the mind of M\=adhava. In support of this surmise, it may be pointed out that Somay\=aj\=\i's methods for computing small arc lengths make use of the M\=adhava-Newton series for the sine function. 

Before taking up N\=\i laka\d n\d tha Somay\=aj\=\i's methods for computing values of the arcsin function, it is illuminating to see how the astronomers and mathematicians of the classical age of Indian astronomy had handled the problem. {\em \=Aryabha\d t\=\i ya} has not considered the problem. 
But Brahmagupta (c. 598--668 CE) and Bh\=askara II (c. 1114--1185 CE) have given approximate expressions for the arcsin function. These expressions are derived from Bh\=askara I's (c. 600--680 CE) approximation function for the sine function. 
Hence, we begin the paper with a discussion on this very interesting approximation function and then consider the associated expression for the arcsin function. 
This is followed by Somay\=aj\=\i's discussion of the methods for computing the arcsin function. There are separate methods for computing  arc lengths corresponding to small jy\=a-s and large jy\=a-s. Both are dealt with in detail. 
Jye\d s\d thadeva's (c. 1500--1575 CE) {\em Yuktibh\=a\d s\=a} contains a novel method for the computation of the circumference of a circle which makes use the arcsin function. We have included a discussion of this also in this chapter. 
\section{Bh\=askara I's approximation to the sine function}\label{BhaskaraApproximationNew}
In this section, we consider a remarkable rational function approximation to the sine function due to Bh\=askara I (c. 600-680). Bhaskara I gave this approximation formula in his {\em Mah\=abh\=askar\=\i ya} and nearly all subsequent Indian mathematicians and astronomers, except those belonging to the Kerala school, have also given equivalent forms of the formula. What is curious and interesting is the fact that that none of the mathematicians of the Kerala school have cared to mention it in their writings; perhaps they thought  they had a more accurate formula in their hands in the form of an infinite series and why bother about an approximation formula!  Neither Bh\=askara I nor any of his successors, true to  their style, have given any rationale for the formula. Historians of mathematics  have come up with several candidate rationales that could be the rationale by which Bh\=askara I originally arrived at the approximation formula. However no candidate has won universal approval!

\subsection{The approximation formula}
The formula is stated in verses 17--19 Chapter VII of {\em Mah\=abh\=askar\=\i ya}. The verses and their English translation are given below (see \cite{Shukla1960} p. 45 for the verses and p. 207 for the English translation; see also \cite{Gupta1967}):
\begin{quote}

\includegraphics[width=6.5cm]{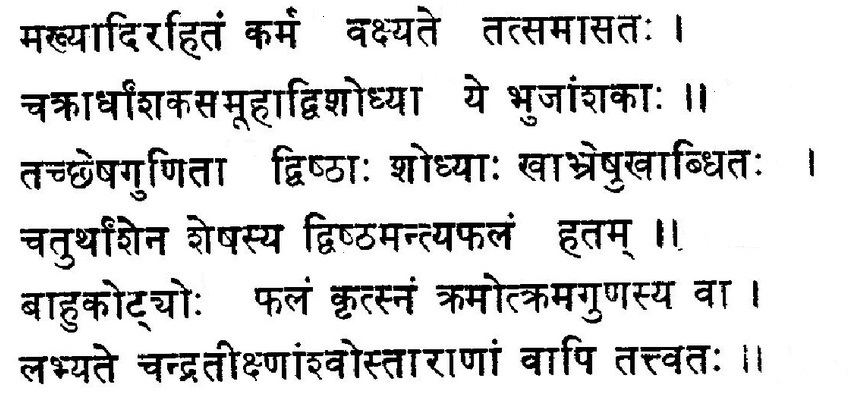}

{\em
makhy\=adirahita\d m karma vak\d syate tatsam\=asata\d h $\scriptstyle \vert$ \\
cakr\=ardh\=a\d m\'sakasam\=uh\=advi\'sodhy\=a ye bhuj\=a\d m\'sak\=a\d h  $\scriptstyle \vert\vert$ \\
tacche\d sugu\d nit\=a dvi\d s\d tha\=a\d h \'sodhy\=a\d h kh\=abhre\d sukh\=abdhita\d ta\d h $\scriptstyle \vert$ \\
caturth\=a\d m\'sena \'se\d sasya dvi\d s\d thamantyaphala\d m hatam $\scriptstyle \vert\vert$\\
b\=ahuko\d tyo\d h phala\d m k\d rtsna\d m kramotkramagu\d nasya v\=a $\scriptstyle \vert$ \\
labhyate candrat\=\i k\d s\d n\=a\d m\'svost\=ar\=a\d n\=a\d m v\=api tattvata\d h $\scriptstyle \vert\vert$   
}

`` I briefly state the rule (for finding the {\em bhujaphala} and {\em kotiphala} etc.) without making use of the Rsine-differences, 225, etc.

Subtract the degrees of the {\em bhuja} (or {\em koti}) from the degrees of half a circle (i.e., $180^\circ$). Then multiply the remainder by the degrees of the {\em bhuja} (or {\em koti}) and put down the result at two places. At one place subtract the result from 40500. By one-fourth of the remainder (thus obtained) divide the result at the other place as multiplied by the {\em antyaphala} (i.e., the epicyclic radius). Thus is obtained the entire {\em bahuphala} (or {\em kotiphala}) for the Sun, Moon, or the star-planets. So also are obtained the direct and inverse Rsines.''
\end{quote}
In the above verses, the reference to ``Rsine-differences, 225, etc.'' is a reference to \=Aryabha\d ta's sine table  for a discussion on this table).
\subsection{Rendering in modern notations}
Let $x^\circ$ be the angle subtended by the arc whose {\em jy\=a}) is sought and let $r$ be the radius of the circle. Bh\=askara I's approximation formula can be expressed in these notations as follows:
\begin{equation}\label{BhaskaraFormulaC}
\text{jy\=a}\,(x^\circ) \approx \frac{r \cdot x(180-x)}{\tfrac{1}{4}[40500 - x(180-x)]}.
\end{equation}
Since $\text{jy\=a}\,(x^\circ) = r \sin(x^\circ)$, we have
\begin{equation}\label{BhaskaraFormula}
\sin(x^\circ) \approx \frac{4x(180-x)}{40500-x(180-x)}.
\end{equation}
This is Bh\=askara I's approximation formula to the sine functions in modern notations. Note that this has to be understood in the following form:
\begin{equation}\label{BhaskaraFormulaA}
\sin(x^\circ) \approx \frac{4\cdot x^\circ\cdot (180^\circ -x^\circ)}{5\cdot 90^\circ\cdot 90^\circ -x^\circ\cdot  (180^\circ -x^\circ)}.
\end{equation}
If angles are measured in radians, Eq.\eqref{BhaskaraFormulaA} can be put in the following form:
\begin{equation}\label{BhaskaraFormulaB}
\sin(\theta) = \frac{16\theta (\pi-\theta)}{5\pi^2 -4\theta(\pi-\theta)}.
\end{equation}
Some properties of the formula can be easily observed.
\begin{enumerate}[\quad (1) ]
\item\label{enum1}
The formula gives the exact values of $\sin (x^\circ)$ for $x=0^\circ, 30^\circ, 90^\circ, 150^\circ$ and $180^\circ$, namely, the values $0, \frac{1}{2}, 1, \frac{1}{2}, 0$.
\item\label{enum2}
The formula is symmetrical about $90^\circ$, that is, the formula gives the same value for $\sin(90^\circ-x^\circ)$ and $\sin(90^\circ + x^\circ)$.
\end{enumerate}
\subsection{Accuracy of the formula}
\begin{figure}
\centering
\includegraphics[width=8cm]{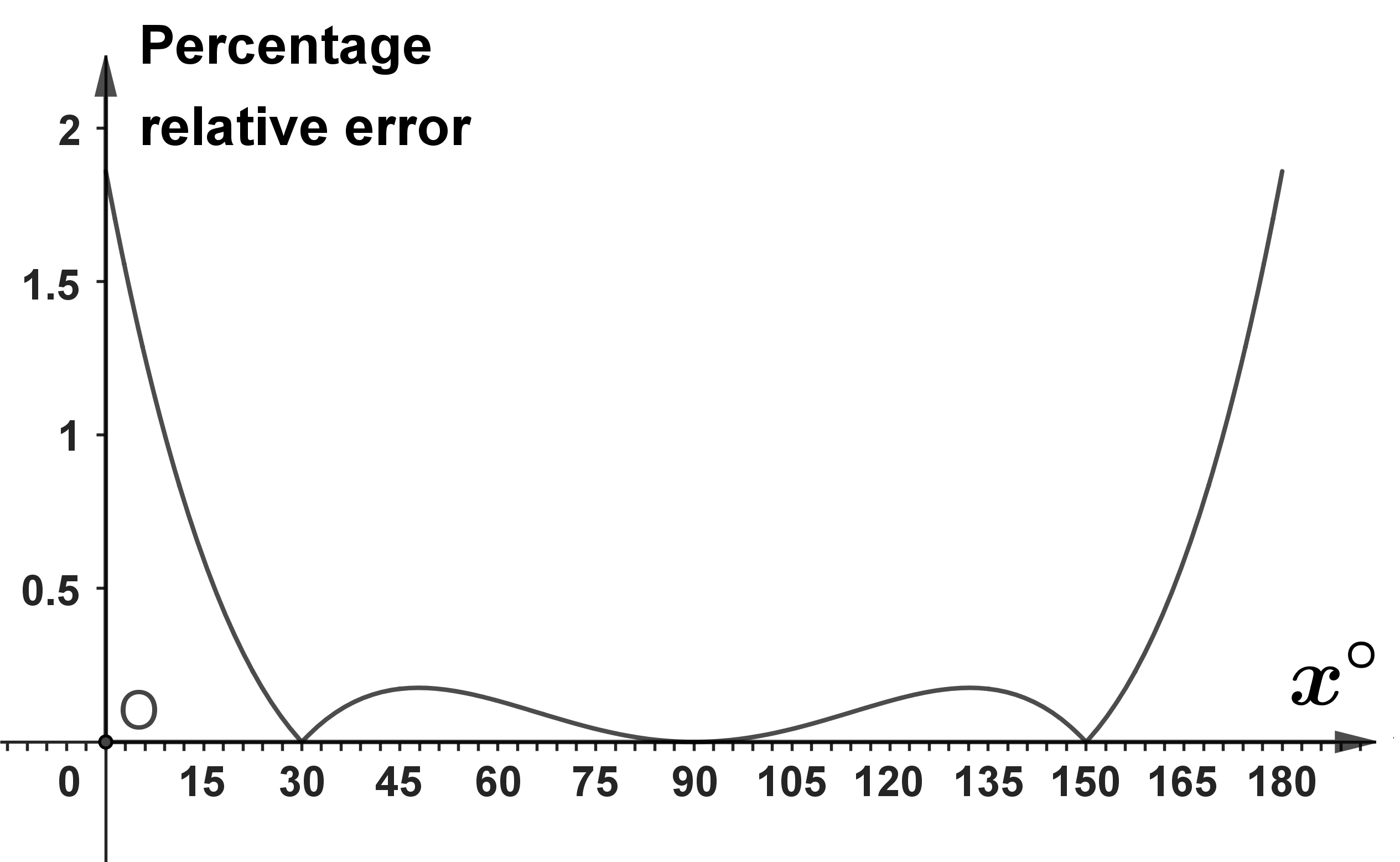}
\caption{Graph of the percentage relative error in Bh\=askara I's approximation formula for the sine function}\label{BhaskaraIerror}
\end{figure}

 Figure \ref{BhaskaraIerror} shows the percentage error in the value of $\sin(x^\circ)$ computed using Bh\=askara I's approximation formula for the sine function. It can be seen that the maximum relative error is less than 1.8\%.
\subsection{On the rationale of the formula}
As already indicated, neither Bh\=askara nor his followers have given the rationale of the sine approximation formula. In the literature on the history of Indian mathematics one can see several modern rationales for the formula (see, for example, \cite{Gupta1967}, \cite{Shirali2011} and \cite{Stroethoff2014}). We present below one such rationale which is an adaption of a rationale given in \cite{Stroethoff2014}.
Our rationale is to look for the simplest function that satisfies the two properties (\ref{enum1}) and (\ref{enum2}) of Bh\=askara I's approximation formula. 
\subsection{The rationale}
The simplest function which satisfies property (\ref{enum2}) and which takes the value $0$ when $x=0, 180$ is 
$$
f(x)=kx(180-x).
$$
Now, $f(30)=\frac{1}{2}$ if $k=\frac{1}{9000}$ and $f(90)=1$ if $k=\frac{1}{8100}$. This suggests that, for $f(x)$ to satisfy property (\ref{enum1}),  $k$ should be of the form $\frac{1}{g(x)}$ with $g(x)$ satisfying property (\ref{enum2}) and with $g(30)=9000$ and $g(90)=8100$. 

The simplest general function which satisfies property (\ref{enum2}) is
$$
g(x)=a+bx(180-x).
$$
Since $g(30)=9000$ and $g(90)=8100$, we must have
\begin{align*}
a+4500b & = 9000,\\
a+8100b & = 8100.
\end{align*}
Solving these equations we get $a=\frac{1}{4}\times 40500$ and $b=-\frac{1}{4}$.

Thus the function
$$
\frac{x(180-x)}{\frac{1}{4}\times 40500 -\frac{1}{4}x(180 - x)}
$$
satisfies the properties (\ref{enum1}) and (\ref{enum2}).
This is precisely Bh\=askara I's sine approximation formula. 
\section{Brahmagupta's formula for arcsin}
Brahmagupta in his {\em Brahmasphu\d tasiddh\=anta} has given a formula for computing the values of the arcsin function. The verses and an English translation thereof are reproduced below (see verses 24--25 Chapter 14 \cite{BrahmaSphuta}; see \cite{Datta1983} for the English translation):
\begin{quote}

\includegraphics[width=5.5cm]{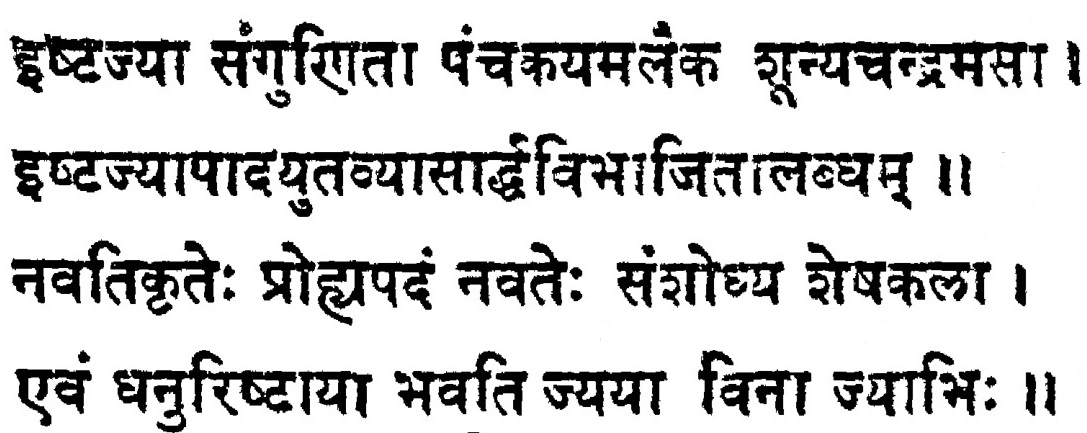}

{\em
i\d s\d tjy\=a sa\d mgu\d nit\=a pa\.mcakayamalaika \'s\=unyacandramas\=a $\scriptstyle\vert$
\\
i\d s\d tajy\=ap\=adayutavy\=as\=arddhavibh\=ajit\=alabdham
$\scriptstyle\vert\vert$
\\
navatik\d rte\d h prohyapada\d m navate\d h sa\d m\'sodhya \'se\d sakal\=a\d h $\scriptstyle\vert$
\\
eva\d m dhanuri\d s\d t\=ay\=a bhavati jyay\=a vin\=a jy\=abhi\d h $\scriptstyle\vert\vert$
}

``Multiply 10125 by the given jy\=a and divide by the quarter of
the given jy\=a plus the radius; subtracting the quotient from the
square of 90, extract the square-root and subtract (the root) from
90; the remainder will be in degrees and minutes; thus will be
found the arc of the given jy\=a without the table of jy\=a-s.''
\end{quote}

\subsection{Rendering in modern notations}
Let $m=\text{jy\=a}\,\,(s)$, where $s$ is an arc, measured in degrees, of a circle of radius $r$, then:
\begin{equation}\label{BrahmaArcSinFormula}
s = 90 - \sqrt{8100- \frac{10125m}{(\frac{m}{4}+r)}}.
\end{equation}
\subsection{Rationale of the formula}
Brahmagupta has not given any rationale for the expression in Eq.\eqref{BrahmaArcSinFormula}. It can be seen that it follows from  from Eq.\eqref{BhaskaraFormulaC}. Rrom Eq.\eqref{BhaskaraFormulaC} we have:
\begin{equation*}
m \approx \frac{r \cdot s(180-s)}{\tfrac{1}{4}[40500 - s(180-s)]}.
\end{equation*}
This can be written as 
$$
s^2 -180s +\frac{10125m}{(\frac{m}{4}+r)}=0.
$$
This is a quadratic equation in $s$ and solving it for $s$ we get precisely the formula in Eq.\eqref{BrahmaArcSinFormula}. The negative  sign of the radical is retained as it is assumed that $s$ is less than $90^\circ$. 
\section{N\=\i laka\d n\d tha Somay\=aj\=\i's  method for computing arcsin of small jy\=a-s}
\subsection{Somay\=aj\=\i's approximate expression for jy\=a\,\,$(s)$}
Somy\=aj\=\i, in {\em Tantras\.ngraha}, has given an 
an approximate expression for jy\=a\,\,$(s)$ in the following verse (see verse 17 Chapter 2 in \cite{Rama2011}):
\begin{quote}
\includegraphics[width=8cm]{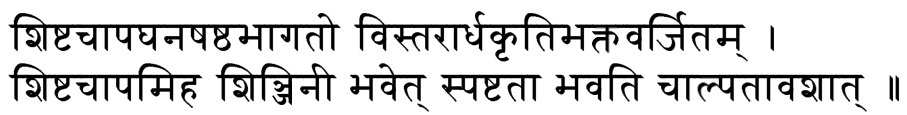}

{\em 
\'si\d s\d tac\=apaghana\d sa\d s\d thabh\=agato vistar\=ardhak\.rtibhaktavarjitam $\scriptstyle\vert$
\\
\'si\d s\d tac\=apamiha \'si\~njin\=\i\  bhavet spa\d s\d tat\=a bhavati c\=alpat\=ava\'s\=at $\scriptstyle\vert\vert$
}

``Divide one-sixth of the cube of the remaining arc by the square of the {\em trijy\=a}. This quantity
when subtracted from the remaining arc becomes the {\em \'si\~njin\=\i\ } (the {\em dorjy\=a} or {\em jy\=a}) corresponding to
the remaining arc). The value is accurate because of the smallness [of the arc].''
\end{quote}
\subsubsection{Rendering in modern notations}
Let $s$ be small arc of a circle of radius $r$, then, the verse states that
\begin{equation}\label{arcsin01}
\text{jy\=a}\,\,(s)=s-\frac{s^3}{6r^2}.
\end{equation}
\subsubsection{Rationale of the result}
Somay\=aj\=\i\ has not given any rationale for the result. However, it can be seen easily that the result follows from the M\=adhava-Newton series for jy\=a\,\,$(s)$, namely, 
\begin{align*}
\text{jy\=a}\,\,(s) = & s - \Big[ s\cdot \frac{s^2}{(2^2+2)r^2} - \Big\{ s\cdot \frac{s^2}{(2^2+2)r^2}\cdot \frac{s^2}{(4^2+4)r^2} - \notag\\
& \phantom{s - }\,\, \Big(s\cdot \frac{s^2}{(2^2+2)r^2}\cdot \frac{s^2}{(4^2+4)r^2}\cdot \frac{s^2}{(6^2+6)r^2} - \Big[ \cdots \Big]\Big)\Big\}\Big],
\end{align*} 
by neglecting terms containing $s^5$ and higher powers of $s$.

\subsection{Formula for computing arcsin function}
Somay\=aj\=\i\ in verse 37 Chapter 2 of {\em Tantrasa\.ngraha} has  indicated that the result in Eq.\eqref{arcsin01} can be used to compute the arc lengths corresponding to a given jy\=a provided that the arcs are small.
However, Somay\=aj\=\i\ did not explicitly state how the value of $s$ is to be computed from this equation. But, fortunately, \'Sankara V\=ariyar  (c. 1500--1560 CE), in his {\em Laghuviv\d rt\=\i} commentary on {\em Tantrasa\.ngraha}, has explained how exactly this can be done (see \cite{Rama2011} pp. 92--93). We rewrite Eq.\eqref{arcsin01} in the following form
\begin{equation}\label{arcsin03}
s=\text{jy\=a}\,\,(s) + \frac{s^3}{6r^2}
\end{equation}
and then, in the right side, we replace $s$ by jy\=a\,\,$(s)$   to get 
\begin{equation}\label{arcsin02}
s\approx \text{jy\=a}\,\,(s) + \frac{(\text{jy\=a}\,\,(s))^3}{6r^2}.
\end{equation}
This is Somay\=aj\=\i's first method for computing values of arc lengths corresponding to small arcs.
\subsection{A polynomial approximation to the arcsin function}
Here we show that Somay\=aj\=\i's expression to compute the arcsin function given by Eq.\eqref{arcsin02} is equivalent to the power series expression for $\arcsin(x)$ truncated at the term containing $x^3$. 

Taking $s=r\theta$ and jy\=a\,\,$(s)=r\sin \theta$, Eq.\eqref{arcsin02} can be written as 
$$
\theta \approx \sin\theta + \frac{\sin^3 \theta}{6}.
$$
Let $x=\sin \theta$ so that $\theta = \arcsin(x)$. Then we have:
\begin{equation}\label{arcsin02a}
\arcsin(x)\approx x+\frac{x^3}{6}.
\end{equation}
Thus Eq.\eqref{arcsin02} is equivalent to the polynomial approximation to $\arcsin(x)$ given by Eq.\eqref{arcsin02a}. Incidentally, the approximation to $\arcsin(x)$ given by Eq.\eqref{arcsin02a} is the approximation to $\arcsin(x)$ obtained by neglecting terms containing $x^5$ and higher powers of $x$ in the power series expression for $\arcsin(x)$, namely,
\begin{equation}\label{arcsin02b}
\arcsin (x)=x+\frac{{{x}^{3}}}{6}+\frac{3 {{x}^{5}}}{40}+\frac{5 {{x}^{7}}}{112}+\frac{35 {{x}^{9}}}{1152}+\cdots.
\end{equation}
\section{\'Sankara V\=ariyar's iterative method for computing approximate values of the arcsin function}
\'Sankara V\=ariyar, in his {\em Laghuviv\d rt\=\i} commentary on {\em Tantras\.ngraha}, after stating the formula given in Eq.\eqref{arcsin02}, has given an iterative method for computing the arc lengths more accurately. The iterative process starts from the approximate relation in Eq.\eqref{arcsin03}, which is rewritten as follows:
\begin{equation}\label{arcsin04}
s=m + \Delta,
\end{equation}
where $m=\text{jy\=a}\,\,(s)$ is known and $s$ is unknown, and where
\begin{equation}\label{arcsin05}
\Delta =\frac{s^3}{6r^2}.
\end{equation}
In Eq.\eqref{arcsin05}, as the first approximation to $s$ we take $s$ as $m$ and compute
$$
\Delta_1=\frac{m^3}{6r^2}
$$
Using this in Eq.\eqref{arcsin04}, we get the first approximation to $s$ as 
$$
s_1=m +\Delta_1.
$$
We next use this value of $s$ in Eq.\eqref{arcsin05} to get the second approximation to $\Delta$ as
$$
\Delta_2=\frac{(m +\Delta_1)^3}{6r^2}.
$$
Using this in Eq.\eqref{arcsin04}, we get the second approximation to $s$ as
$$
s_2=m +\Delta_2.
$$
The process is continued until two successive approximations to $s$ are nearly equal. 

\subsection{\'Sankara V\=ariyar's Algorithm to compute arc lengths}\label{SankaraVariyarAlgorithm}
The above described procedure can be summarised as an algorithm thus. Let $m$ be the jy\=a of an arc $s$ of a circle of radius $r$. Given $m$, to compute a sequence of approximations $s_i$ to $s$ and to generate an approximate value of $s$:
\begin{enumerate}[\qquad{Step} 1.]
\item
Set $s_0=m$, $\Delta_0=0$.
\item
For $i=1,2,3, \ldots$, repeat the following until successive values of $s_i$-s are nearly equal:
\begin{enumerate}
\item
$\Delta_i= \frac{(m+\Delta_{i-1})^3}{6r^2}.$
\item
$s_i=m+\Delta_i.$
\end{enumerate}
\item
When $s_i\approx s_{i+1}$, then $s\approx s_i$.
\end{enumerate}
\subsubsection*{Remarks}
The above algorithm to compute arc lengths does not generate successively better approximations to the true value of the arc length. Moreover, the sequence $\{s_i\}$ does not converge to the true arc length $s$. In fact we have:
\begin{align*}
s_0& = m,\\
s_1&= m+\frac{m^3}{6r^2},\\
s_2&=m+  \frac{{{m}^{3}}}{6 {{r}^{2}}}+\frac{{{m}^{5}}}{12 {{r}^{4}}}+\frac{{{m}^{7}}}{72 {{r}^{6}}}+\frac{{{m}^{9}}}{1296 {{r}^{8}}},\\
s_3&=m+ \frac{{{m}^{3}}}{6 {{r}^{2}}}+\frac{{{m}^{5}}}{12 {{r}^{4}}}+\frac{{{m}^{7}}}{18 {{r}^{6}}}+\frac{7 {{m}^{9}}}{324 {{r}^{8}}}+ \text{ higher powers of $m$},\\
s_4&=m+\frac{{{m}^{3}}}{6 {{r}^{2}}}+\frac{{{m}^{5}}}{12 {{r}^{4}}}+\frac{{{m}^{7}}}{18 {{r}^{6}}}+\frac{55 {{m}^{9}}}{1296 {{r}^{8}}}+\text{ higher powers of $m$}.
\end{align*}
These may be compared with the power series expansion of the arcsin function given in Eq.\eqref{arcsin02b}.
\subsection{Illustrative examples}
\subsubsection{Illustrative example 1}
Let us compute the arc whose jy\=a is $m=224^\prime\,\,50^{\prime\prime}\,\,22^{\prime\prime\prime}= 809422^{\prime\prime\prime}$. This is the first entry in M\=adhava's sine table and it corresponds to jy\=a\,\,$(225^\prime)$. We take the radius of the circle as $r=3437^\prime\,\,44^{\prime\prime}\,\,48^{\prime\prime\prime}= 12375888^{\prime\prime\prime}$. 
\begin{enumerate}[\qquad{Step} 1.]
\item
As initial approximations, we take $s_0=809422^{\prime\prime\prime}$ and $\Delta_0=0$.
\item
Results of computations are taken as rounded to the nearest integer.

Iteration 1:

$\Delta_1=\frac{(m+\Delta_0)^3}{6r^2}=\frac{(809422^{\prime\prime\prime}+0)^3}{6\times(12375888^{\prime\prime\prime})^2}= 577^{\prime\prime\prime}$

$s_1=  m +\Delta_1 = 809422^{\prime\prime\prime} + 577^{\prime\prime\prime} = 809999^{\prime\prime\prime}$

Iteration 2:

$\Delta_2=\frac{(m+\Delta_1)^3}{6r^2}=\frac{(809422^{\prime\prime\prime} + 577^{\prime\prime\prime})^3}{6\times(12375888^{\prime\prime\prime})^2}=578^{\prime\prime\prime}$

$s_2 = m+\Delta_2 = 809422^{\prime\prime\prime} + 578^{\prime\prime\prime} = 810000^{\prime\prime\prime}$

Iteration 3:

$\Delta_3=\frac{(m+\Delta_2)^3}{6r^2}=\frac{(809422^{\prime\prime\prime} + 578^{\prime\prime\prime})^3}{6\times(12375888^{\prime\prime\prime})^2}=578^{\prime\prime\prime}$

$s_3 = m+\Delta_3 = 809422^{\prime\prime\prime} + 578^{\prime\prime\prime} = 810000^{\prime\prime\prime}$

\item
$s_2=s_3= 810000^{\prime\prime\prime}$ and so $s=810000^{\prime\prime\prime} = 225^\prime$. This agrees with the true value of $s$.
\end{enumerate}

\subsubsection{Illustrative example 2}
Let us compute the arc whose jy\=a is $m=448^\prime\,\,42^{\prime\prime}\,\,58^{\prime\prime\prime}= 1615378^{\prime\prime\prime}$. This is the second entry in M\=adhava's sine table and it corresponds to jy\=a\,\,$(450^\prime)$.
\begin{enumerate}[\qquad{Step} 1.]
\item
As initial approximations, we take $s_0= 1615378^{\prime\prime\prime}$ and $\Delta_0=0$.
\item
Results of computations are taken as rounded to the nearest integer.

Iteration 1:

$\Delta_1=\frac{(m+\Delta_0)^3}{6r^2}= 4587^{\prime\prime\prime}$

$s_1=  m +\Delta_1  = 1619965^{\prime\prime\prime}$

Iteration 2:

$\Delta_2=\frac{(m+\Delta_1)^3}{6r^2}=4626^{\prime\prime\prime}$

$s_2 = m+\Delta_2  = 1620004^{\prime\prime\prime}$

Iteration 3:

$\Delta_3=\frac{(m+\Delta_2)^3}{6r^2}=4626^{\prime\prime\prime}$

$s_3 = m+\Delta_3 = 1620004^{\prime\prime\prime}$

\item
$s_2=s_3= 1620004^{\prime\prime\prime}$ and so $s=1620004^{\prime\prime\prime} = 450^\prime\,\,0^{\prime\prime}\,\,4^{\prime\prime\prime}$. In this case, there is a difference of $4^{\prime\prime\prime}$ from the true value of $s$. 
\end{enumerate}
\subsection{\'Sankara V\=ariyar's algorithm throws up a core integer sequence!}
The algorithm described in Section \ref{SankaraVariyarAlgorithm} has completely unexpectedly thrown up an important integer sequence which appears in a large number of combinatorial problems. The Online Encyclopedia of Integer Sequences has characterised this as a ``core integer sequence''. To see how the particular integer sequence emerges via \'Sankara V\=ariyar's algorithm, let us reformulate the algorithm in the following form: we write $m=x$, $t=\frac{1}{6r^2}$, $\Delta_i=y_i$, $s_0=x$ and $y_0=0$. Then we have, for $i=1,2,\ldots$:
\begin{align*}
y_{i}&=t(x+y_{i-1})^3\\
s_{i}&=x+y_{i-1}.
\end{align*}
Let us compute a few values of $s_i$:
\begin{align*}
s_0 & =
(x), \\
s_1 & = 
(x+t {{x}^{3}}),
\\
s_2 & = 
(x+t {{x}^{3}}+3 {{t}^{2}} {{x}^{5}})\;+3 {{t}^{3}} {{x}^{7}}+{{t}^{4}} {{x}^{9}},
\\
s_3 & =
(x+t {{x}^{3}}+3 {{t}^{2}} {{x}^{5}}+12 {{t}^{3}} {{x}^{7}})\;+28 {{t}^{4}} {{x}^{9}}+57 {{t}^{5}} {{x}^{11}}+96 {{t}^{6}} {{x}^{13}}+ \text{ higher powers of $x$},
\\
s_4 & = 
(x+t {{x}^{3}}+3 {{t}^{2}} {{x}^{5}}+12 {{t}^{3}} {{x}^{7}}+55 {{t}^{4}} {{x}^{9}})\;+192 {{t}^{5}} {{x}^{11}}+618 {{t}^{6}} {{x}^{13}}+\text{ higher powers of $x$},
\\
s_5 & = 
(x+t {{x}^{3}}+3 {{t}^{2}} {{x}^{5}}+12 {{t}^{3}} {{x}^{7}}+55 {{t}^{4}} {{x}^{9}}+273 {{t}^{5}} {{x}^{11}})\;+1185 {{t}^{6}} {{x}^{13}}+\text{ higher powers of $x$},
\\
s_6 & =
(x+t {{x}^{3}}+3 {{t}^{2}} {{x}^{5}}+12 {{t}^{3}} {{x}^{7}}+55 {{t}^{4}} {{x}^{9}}+273 {{t}^{5}} {{x}^{11}}+1428 {{t}^{6}} {{x}^{13}})\;+\text{ higher powers of $x$}.
\end{align*}
Consider the integers that appear as coefficients in the parenthesised expressions in the $s_i$-s:
\[
1, 1, 3, 12, 55, 273, 1428, \ldots.
\]
This is precisely the sequence A001764 in the Online Encyclopedia of Integer Sequences (OIES) (see \cite{OEISA001764}). This sequence, which appears in a large number of combinatorial problems, has been characterised as a core sequence in the Encyclopedia.  As per the Encyclopedia, the $j$-th entry $a_j$ in the above integer sequence is given by 
\[
a_j=\frac{(3j)!}{j!(2j+1)!}.
\]
Thus we have:
\[
s_n= (a_0x+a_1tx^3+a_2t^2x^5+\cdots + a_nt^nx^{2n+1})\;+\text{ higher powers of $x$}.
\]
Moreover, there is a closed form expression for the infinite series $\sum_{n=0}^\infty a_nt^nx^{2n+1}$ (for example, see \cite{TomCopeland}), namely,
\[
\frac{2}{\sqrt{3 t}} \sin{\Big( \frac{1}{3} \operatorname{arcsin}\Big( \frac{3\sqrt{3 t}}{2}x\Big) \Big) }
=
x+t {{x}^{3}}+3 {{t}^{2}} {{x}^{5}}+12 {{t}^{3}} {{x}^{7}}+55 {{t}^{4}} {{x}^{9}}+273 {{t}^{5}} {{x}^{11}}+\cdots .
\]

\section{Computing arcsin function using a lookup table}
N\=\i laka\d n\d tha Somay\=aj\=\i\ has suggested that a lookup table can be used to find the arc lengths corresponding to certain pre-defined jy\=a values. This method has a severe limitation: it cannot be used to find the arc lengths corresponding to any arbitrary jy\=a values. The table lists only 24 jy\=a values the minimum being $105^\prime\,\, 43^{\prime\prime}$ and the maximum being $304^\prime\,\,36^{\prime\prime}$. The difference between a pair successive jy\=a values is not a constant. 

\subsection{The rationale of the method}
To begin with, we assume that the arcs and the corresponding jy\=a-s are all small. Let $s$ be a small arc and let jy\=a\,\,$(s)=m$. We are given $m$ and we have to find $s$. From Eq.\eqref{arcsin03}, we have
\begin{equation}\label{arcsin06}
s-m\approx \frac{m^3}{6r^2}.
\end{equation}
Somay\=aj\=\i\ assumes that $m$ and $r$ are in minutes so that $s-m$ given by Eq.\eqref{arcsin06} is also in minutes. If 
\begin{equation}\label{arcsin07}
s-m=k^{\prime\prime}=(\tfrac{k}{60})^\prime,
\end{equation}
 then we have
\begin{equation}
\frac{k}{60}=\frac{m^3}{6r^2}.
\end{equation}
From this, we get
\begin{equation}\label{arcsin09a}
m=\Big(\frac{kr^2}{10}\Big)^{\frac{1}{3}}\,\,\text{minutes}.
\end{equation}
Using this in Eq.\eqref{arcsin07}, we get
\begin{equation}\label{arcsin08}
s= \Big[\frac{k}{60} + \Big(\frac{kr^2}{10}\Big)^{\frac{1}{3}}\Big] \,\,\text{minutes}.
\end{equation}
This means that, if the difference between arc length and jy\=a is $k^{\prime\prime}$, then the arc length is given Eq.\eqref{arcsin08}. The lookup table gives the 
the values of $s$  for $k=1,2,\ldots,24$.
\subsection{The table}
Table \ref{arcsin_lookup_table} lists the values of the arc lengths for various values of the difference between the arc lengths and the corresponding jy\=a-s as given in the {\em Laghuviv\d rt\=\i} commentary of {\em Tantrasa\.ngraha} (see \cite{Rama2011} pp. 94--95). The {\em Laghuviv\d rt\=\i} commentary lists the values of $m=\text{jy\=a}\,\,(s)$ given by Eq.\eqref{arcsin09a} for $k=1,2,\ldots,24$ whereas Table \ref{arcsin_lookup_table} lists the values of $\text{jy\=a}\,\,(s)$ as well as the values of $s$ given by Eq.\eqref{arcsin08}. In the table, we have also given the {\em ka\d tapay\=adi} encodings of the values of jy\=a\,\,$(s)$-s. These encodings are taken from the following  verses given in the {\em Laghuviv\d rt\=\i} commentary.
\begin{quote}
\includegraphics[width=6.5cm]{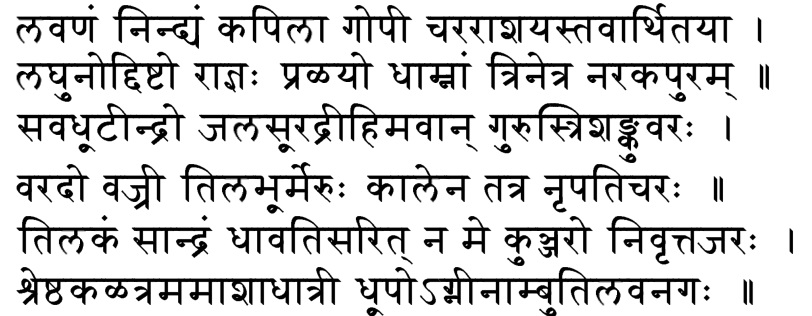}

{\em
lava\d na\d m nindya\d m kapil\=a gop\=\i\ carar\=a\'sayastav\=arthitay\=a $\scriptstyle\vert$
\\
laghunoddi\d s\d to r\=aj\~na\d h pra\d layo dh\=amn\=a\d m trinetra narakapuram $\scriptstyle\vert\vert$
\\
savadh\=u\d t\=\i ndro jalas\=uradr\=\i himav\=an gurustri\'sa\.nkuvara\d h $\scriptstyle\vert$
\\
varado vajr\=\i\ tilabh\=urmeru\d h k\=alena tatra n\d rpaticara\d h $\scriptstyle\vert\vert$
\\
tilaka\d m s\=andra\d m dh\=avatisarit na me ku\~njaro niv\d rttajara\d h $\scriptstyle\vert$
\\
\'sre\d s\d thaka\d latramam\=a\'s\=adh\=atr\=\i\ dh\=upo'gn\=\i n\=ambutilavanaga\d h $\scriptstyle\vert\vert$
}
\end{quote}
\renewcommand{\arraystretch}{1.1}
\begin{table}[ht]
\footnotesize
\centering
\begin{tabular}{c|c|cc|cc}
\hline
$k$ & jy\=a\,$(s)$ & \multicolumn{2}{c|}{jy\=a\,$(s)$} & \multicolumn{2}{c}{$s$} \\
\hline
  $\prime\prime$ & (in {\em ka\d tapay\=adi}) & $\prime$ & $\prime\prime$ & $\prime$ & $ \prime\prime$ 
\\
\hline
1  & {\em lava\d na\d m nindya\d m} & 105 & 43 & 105 & 44 \\
2  & {\em kapil\=a gop\=\i\ }       & 133 & 11 & 133 & 13 \\
3  & {\em carar\=a\'saya }          & 152 & 26 & 152 & 29 \\
4  & {\em stav\=arthitay\=a}        & 167 & 46 & 167 & 50 \\
5  & {\em laghunoddi\d s\d to }     & 180 & 43 & 180 & 48 \\
6  & {\em r\=aj\~na\d h pra\d layo} & 192 & 02 & 192 & 08 \\
7  & {\em dh\=amn\=a\d m trinetra } & 202 & 08 & 202 & 15 \\
8  & {\em narakapuram}              & 211 & 20 & 211 & 28 \\
9  & {\em savadh\=u\d t\=\i ndro}   & 219 & 47 & 219 & 56 \\
10 & {\em jalas\=uradr\=\i }        & 227 & 38 & 227 & 48 \\
11 & {\em himav\=an guru }          & 234 & 58 & 235 & 09 \\
12 & {\em stri\'sa\.nkuvara\d h}    & 241 & 52 & 242 & 04 \\
13 & {\em varado vajr\=\i\ }        & 248 & 24 & 248 & 37 \\
14 & {\em tilabh\=urmeru\d h}       & 254 & 36 & 254 & 50 \\
15 & {\em k\=alena tatra }          & 260 & 31 & 260 & 46 \\
16 & {\em n\d rpaticara\d h }       & 266 & 10 & 266 & 26 \\
17 & {\em tilaka\d m s\=andra\d m } & 271 & 36 & 271 & 53 \\
18 & {\em dh\=avatisarit  }         & 276 & 48 & 277 & 06 \\
19 & {\em na me ku\~njaro }         & 281 & 50 & 282 & 09 \\
20 & {\em niv\d rttajara\d h}       & 286 & 40 & 287 & 00 \\
21 & {\em \'sre\d s\d thaka\d latra}& 291 & 22 & 291 & 43 \\
22 & {\em mam\=a\'s\=adh\=atr\=\i\ }& 295 & 55 & 296 & 17 \\
23 & {\em dh\=upo'gn\=\i n\=a }     & 300 & 18 & 300 & 41 \\
24 & {\em mbutilavanaga\d h}        & 304 & 36 & 305 & 00 \\
\hline
\end{tabular}
\caption{Lookup table for computing values of the arcsin function}\label{arcsin_lookup_table}
\end{table}
\renewcommand{\arraystretch}{1}

The values given by the {\em Laghuviv\d rt\=\i} commentary  are remarkably accurate. For example, for $k=24$, the commentary gives $s=305^\prime\,\,0^{\prime\prime}$ whereas a modern computation using {\em Maxima} software, with $r=\frac{21600}{2\pi}$, yielded $s=304^\prime\,\,58.03^{\prime\prime}$, a difference of nearly $2^{\prime\prime}$ only.
\subsection{How to use the lookup table}
The lookup table gives the values of the arc lengths corresponding to a predefined set of jy\=a values. If the jy\=a whose arc length is to be determined happens to be very close to one of the values listed in the table,
then the corresponding arc length can be read off from the table. For example, if $\text{jy\=a}\,\,(s)=200^\prime$ which is very close to the value $202^\prime\,\, 08^{\prime\prime}$ listed in the table, then we may take $s$ as the corresponding arc length  given in the table,  namely, $202^\prime\,\,15^{\prime\prime}$. 

It is interesting to speculate why the table contains only values up to $k=24^{\prime\prime}$. One important application of the table might have been the computation of the true longitudes of the Moon which  involved the computation of an expression of the form (see \cite{Rama2011} p. 90)
$$
\arcsin\big(\tfrac{7}{80}\sin(\theta_0 -\theta_m)\big).
$$
The maximum value of $\tfrac{7}{80}\sin(\theta_0 -\theta_m)$ is
$$
\tfrac{7}{80}\times (3437^\prime\,\, 44^{\prime\prime}\,\,48^{\prime\prime\prime})=300^\prime\,\, 48^{\prime\prime}\,\,10^{\prime\prime\prime}.
$$
Hence in the computation of the true longitudes of the Moon, one would be required to compute arc lengths corresponding to jy\=a-s less than this maximum value only. Note that the value of jy\=a corresponding to $k=23^{\prime\prime}$ is $300^\prime\,\, 18^{\prime\prime}$ and that corresponding to $k=24^{\prime\prime}$ is $304^\prime \,\, 36^{\prime\prime}$, and that the maximum value lies between these two values.

\section{N\=\i laka\d n\d tha Somay\=aj\=\i's  method for computing arcsin of large jy\=a-s }
N\=\i laka\d n\d tha Somay\=aj\=\i\ in his {\em tantrasa\.ngraha} has given an approximation formula for the difference between two arcs in terms of the jy\=a-s and kojy\=a-s of the arcs. This formula is then used to compute the value of the arcsin function corresponding to a large value of jy\=a. 

\subsection{Somay\=aj\=\i's formula}
N\=\i laka\d n\d tha Somay\=aj\=\i's formula appears in verse 14 Chapter 2 of {\em tantrasa\.ngraha} (see \cite{Rama2011} p. 68). 
\begin{quote}

\includegraphics[width=6cm]{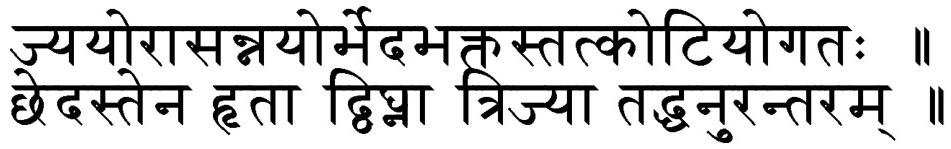}

{\em 
jyayor\=asannayorbhedabhaktastatko\d tiyogata\d h
 $\scriptstyle\vert$
 \\
chedastena h\d rt\=advighn\=a trijy\=a taddhanurantaram
$\scriptstyle\vert\vert$
}

``The sum of the kojy\=a-s divided by the difference of those two jy\=a-s, which are close to each
other, forms the {\em cheda} (divisor). Twice the {\em trijy\=a} divided by this is the difference between
the corresponding arcs.''
\end{quote}

\subsubsection{Rendering in modern notations}
Let $s_1$ and $s_2$ be two arcs, with $s_2>s_1$, of a circle of radius $r$. Then the quoted verse says:
\begin{equation}\label{arcsin09}
s_2-s_1 = \frac{2r}{\left( \dfrac{\text{kojy\=a}\,\,(s_2) + \text{kojy\=a}\,\,(s_1)}{\text{jy\=a}\,\,(s_2) - \text{jy\=a}\,\,(s_1)}\right)}.
\end{equation}
\subsubsection{Equivalent mathematical result}
Let the arcs $s_1$ and $s_2$ subtend angles $\theta_1$ and $\theta_2$ at the centre of the circle. Then:
\begin{gather*}
 \text{jy\=a}\,\,(s_1)=r\sin\theta_1,\quad \text{jy\=a}\,\,(s_2)=r\sin\theta_2,\\
  \text{kojy\=a}\,\,(s_1)=r\cos\theta_1,\quad \text{kojy\=a}\,\,(s_2)=r\cos\theta_2.
\end{gather*}
Using these in Eq.\eqref{arcsin09}, we have:
\begin{equation}\label{arcsin10a}
s_2-s_1 = \frac{2r}{\left(\dfrac{\cos\theta_2+\cos\theta_1}{\sin\theta_2-\sin\theta_1}\right)}.
\end{equation}
\subsubsection{Remarks}
The mathematical basis of Eq.\eqref{arcsin09} is a well-known result of elementary trigonometry which states that 
\begin{equation}\label{arcsin11}
\tan \theta \approx \theta \text{ when $\theta$ is small}.
\end{equation}
This can be seen by expressing the right side of Eq.\eqref{arcsin10a} in the following form:
\begin{align}
\frac{2r}{\left(\dfrac{\cos\theta_2+\cos\theta_1}{\sin\theta_2-\sin\theta_1}\right)}
& = 2r \frac{(\sin\theta_2-\sin\theta_1)}{\cos\theta_2+\cos\theta_1} \notag \\
& = 2r \dfrac{2\cos\frac{\theta_2+\theta_1}{2}\sin\frac{\theta_2-\theta_1}{2}}{2 \cos\frac{\theta_2+\theta_1}{2}\cos\frac{\theta_2-\theta_1}{2}}\notag \\
& = 2r \tan \left(\frac{\theta_2-\theta_1}{2}\right)\notag \\
&=2r \tan \left(\frac{s_2-s_1}{2r}\right).\label{arcsin13}
\end{align}
Now, using Eq.\eqref{arcsin11}, we have:
\[
2r\tan\left(\frac{s_2-s_1}{2r}\right)\approx 2r\Big(\frac{s_2-s_1}{2r}\Big)=s_2-s_1.
\]
Thus Eq.\eqref{arcsin10a}, and consequently Eq.\eqref{arcsin09} also,  follows from Eq.\eqref{arcsin11}.
\subsection{Rationale of the formula}
{\em Tantrasa\.ngraha} does not give any rationale for 
Eq.\eqref{arcsin09}. The editors of a modern edition of 
{\em Tantrasa\.ngraha} have suggested that a geometrical argument could be the method by which Indian astronomers arrived at the result (see \cite{Rama2011} pp. 69--70). We present such an argument below which is a slight variant of the argument given in \cite{Rama2011}. 
\begin{figure}[ht]
\centering
\includegraphics[width=6cm]{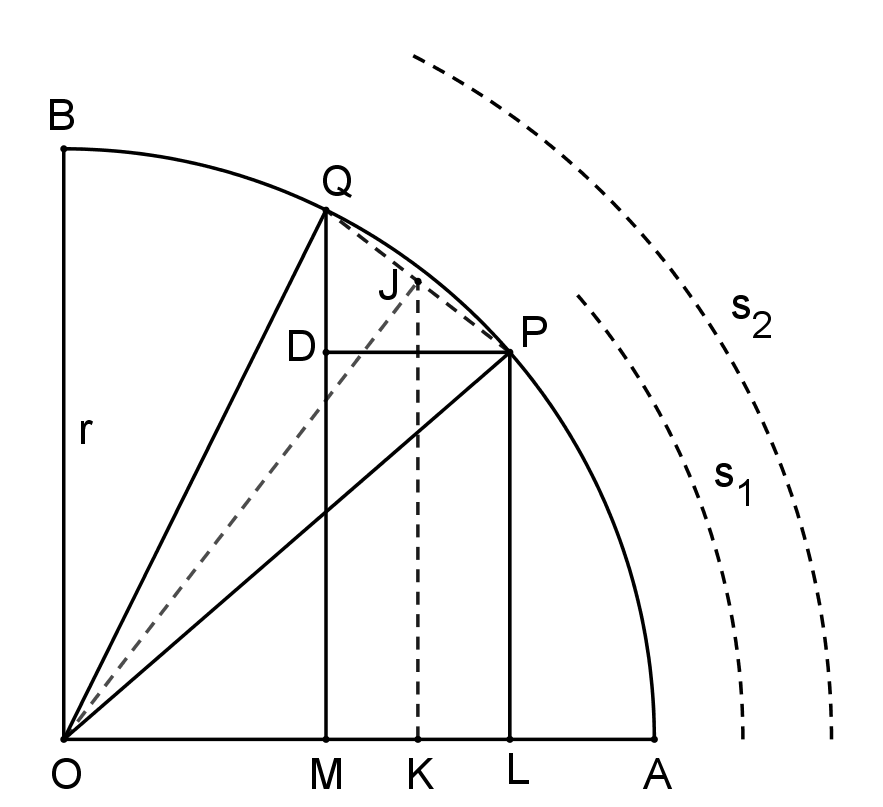}
\caption{Diagram to derive the formula for computing arc lengths}\label{ArcLength02}
\end{figure}

In Figure \ref{ArcLength02}, we assume that $s_2-s_1 =\arc{PQ}$ is small and so we take 
\begin{equation}\label{arcsin13a}
s_2-s_1=\arc{PQ}=PQ.
\end{equation}
We shall now compute the length of the chord $PQ$. 

Let $J$ be the midpoint of the chord $PQ$. We drop perpendiculars $PL, QM, JK$ to $OA$ and $PD$ to $QM$. $K$ is the midpoint of $ML$. From the figure we have:
\begin{align*}
QD & = QM-PL\\
   & = \text{jy\=a}\,\,(s_2) - \text{jy\=a}\,\,(s_1),\\
OK & = \tfrac{1}{2}(OL+OM)\\
   & = \tfrac{1}{2}( \text{kojy\=a}\,\,(s_2)+\text{kojy\=a}\,\,(s_1)) .  
\end{align*}
Since $J$ is the midpoint of $PQ$, $OJ$ is perpendicular to $PQ$. It follows that the triangles $PQD$ and $JOK$ are similar. hence:
$$ 
\frac{PQ}{JO}=\frac{QD}{OK}.
$$
Using the expressions for $QD$ and $OK$, we have:
$$
PQ = \frac{2\cdot JO}{ \left(\dfrac{\text{kojy\=a}\,\,(s_2) + \text{kojy\=a}\,\,(s_1)}{\text{jy\=a}\,\,(s_2) - \text{jy\=a}\,\,(s_1)}\right)}.
$$
Now, $JO\approx r$ and so we have
\begin{equation}\label{arcsin14}
PQ \approx \frac{2r}{ \left(\dfrac{\text{kojy\=a}\,\,(s_2) + \text{kojy\=a}\,\,(s_1)}{\text{jy\=a}\,\,(s_2) - \text{jy\=a}\,\,(s_1)}\right)}.
\end{equation}
Combining Eq.\eqref{arcsin13a} and Eq.\eqref{arcsin14} we get
$$
s_2-s_1
\approx \frac{2r}{ \left(\dfrac{\text{kojy\=a}\,\,(s_2) + \text{kojy\=a}\,\,(s_1)}{\text{jy\=a}\,\,(s_2) - \text{jy\=a}\,\,(s_1)}\right)}.
$$
\subsection{Computation of arcsin of large jy\=a-s}
Let $s$ be an unknown arc whose jy\=a is known, say, jy\=a\,\,$(s)=m$. Using the value of jy\=a\,\,$(s)$, the value of kojy\=a\,\,$(s)$ can be calculated using the the following result:
\[
\text{kojy\=a}\,\,(s)=\sqrt{r^2-(\text{jy\=a}\,\,(s))^2}=\sqrt{r^2-m^2}.
\]
Now, let $m$ be between two consecutive jy\=a values, say jy\=a\,\,$(s_1)$ and jy\=a\,\,$(s_2)$, listed in a table of jy\=a-s, say M\=adhava's sine table,  and let $s_1<s_2$. Let $m$ be closer to jy\=a\,\,$(s_1)$. Using Eq.\eqref{arcsin09}, we have:
\begin{align*}
s-s_1 & =\frac{2r}{\left(\dfrac{\text{\text{kojy\=a}\,\,(s) + kojy\=a}\,\,(s_1)}{\text{jy\=a}\,\,(s)-\text{jy\=a}\,\,(s_1)}\right)}\\
&= \frac{2r}{\left(\dfrac{\sqrt{r^2-m^2} + \text{kojy\=a}\,\,(s_1)}{m-\text{jy\=a}\,\,(s_1)}\right)}\\
&= p,\text{ say},
\end{align*}
from which we get $s=s_1+p$. If $m$ is closer to jy\=a\,\,$(s_2)$, we have to calculate $s_2-s$ and subtract it from $s_2$ to get $s$.
\subsubsection{Error estimate}
Let the true value of $s-s_1$ be $p^\ast$. Then, by Eq.\eqref{arcsin13}, the approximate value of $s-s_1$ is 
$$
p=2r\tan \big(\tfrac{p^\ast}{2r}\big).
$$
The error in the approximate value is 
$$
|p-p^\ast| = \big| 2r\tan \big(\tfrac{p^\ast}{2r}\big)- p^\ast \big|.
$$
If we are using M\=adhava's sine table, the maximum value of $p^\ast$ is $225^\prime$. So, taking $r=3437^\prime\,\,44^{\prime\prime}\,\,48^{\prime\prime\prime}$, the maximum error in the computed value of $p^\ast$ is 
$$
\big|2r\tan \big(\tfrac{225^\prime}{2r}\big) - 225^\prime \big| = 4^{\prime\prime}\,\, 49^{\prime\prime\prime}.
$$
\subsection{Illustrative example}
To illustrate the method with a numerical example, let us calculate arcsin\,\,$(3000^\prime) = s$. So, here $m=\text{jy\=a}\,\,(s)=3000^\prime$. We perform all computations in the units of arc-thirds and hence the radius $r$ will be taken as 
$$
r=(3437\cdot 60\cdot 60 +44\cdot 60 +48)^{\prime\prime\prime}=12375888^{\prime\prime\prime}
$$
and $m$
as
$$
m=(3000\cdot 60\cdot 60)^{\prime\prime\prime}=10800000^{\prime\prime\prime}
$$
We first calculate koj\=a\,\,$(s)$:
\begin{align*}
\text{koj\=a}\,\,(s)
& = \sqrt{r^2-m^2}\\
& = \sqrt{12375888^2 - 10800000^2}\\
& = \sqrt{36522603788544}\\
& = 6043393^{\prime\prime\prime}\text{ rounded to nearest integer}.
\end{align*}
In M\=adhava's sine table, the jy\=a value $3000^\prime$ lies between $\text{jy\=a}\,\,(3600^\prime)=2977^\prime\,\,10^{\prime\prime}\,\,34^{\prime\prime\prime}$  and $\text{jy\=a}\,\,(2825^\prime)=3083^\prime\,\,13^{\prime\prime}\,\,17^{\prime\prime\prime}$ and is closer to $\text{jy\=a}\,\,(3600^\prime)$. Hence we take $s_1=3600^\prime$ and 
\[
\text{jy\=a}\,\,(s_1)=(2977\cdot 60\cdot 60+10\cdot 60 + 34)^{\prime\prime\prime} = 10717834^{\prime\prime\prime}.
\]
Further, we have:
\begin{align*}
\text{kojy\=a}\,\,(s_1)
& = \text{kojy\=a}\,\,(3600^\prime)\\
& = \text{jy\=a}\,\,(5400^\prime - 3600^\prime)\\
& = \text{jy\=a}\,\,(1800^\prime)\\
& = (1718\cdot 60\cdot 60 + 52\cdot 60 + 24)^{\prime\prime\prime}\\
& = 6187944^{\prime\prime\prime}.
\end{align*}
Therefore, we have
\begin{align*}
p
& = s-s_1 \\
&= \frac{2\cdot 12375888}{\left(\frac{6043393 + 6187944 }{9000000 - 10717834}\right)}\\
& = 166274^{\prime\prime\prime},
\end{align*}
and hence
\begin{align*}
s 
& =s_1+p\\
& = 3600^\prime + 166274^{\prime\prime\prime}\\
& = 3646^\prime\,\, 11^{\prime\prime}\,\,14^{\prime\prime\prime}.
\end{align*}
This is amazingly accurate because by a modern computation we see that 
\begin{align*}
\text{jy\=a}\,\,(3646^\prime\,\, 11^{\prime\prime}\,\,14^{\prime\prime\prime}) 
& = r \sin( 3646^\prime\,\, 11^{\prime\prime}\,\,14^{\prime\prime\prime})\\
& = 12375888^{\prime\prime\prime}  \sin( 3646^\prime\,\, 11^{\prime\prime}\,\,14^{\prime\prime\prime})\\
& = 10800001^{\prime\prime\prime}
\end{align*}
 and the difference of this  from the value of $m$ is only $1^{\prime\prime\prime}$.
\section{A novel method for computing the circumference of a circle}
Two different  methods, both due to M\=adhava, for the computation of the circumference of a given circle are well known. One of the methods  makes use of the several infinite series expressions for $\pi$ and the correction terms and the other  is an iterative procedure using geometrical ideas. In this section we present a third method for the computation of the circumference. Even though the present method produces only an approximate value of the circumference and is not of much practical significance, we have taken it for presentation for the beauty and elegance of the mathematics involved. (Does this indicate that the astronomers and mathematicians of the Kerala school enjoyed ``doing mathematics'' for its own sake?) The method and its rationale appears in {\em Yuktibh\=a\d s\=a} (see {\em Ga\d nita-Yukti-Bh\=a\d s\=a} Section 7.6; \cite{Sarma2008} pp. 233-234) and is presented there as an application of the M\=adhava-Newton series for jy\=a.  It is also an application of the {\em j\=\i ve-paraspara ny\=aya}-s due to M\=adhava, the addition and subtraction rules for the sine and cosine functions. Though we are sure that this method is a product of the Kerala school of mathematics, we are not sure who was the original inventor of the method.
\subsection{Computation of the circumference of a circle}
\begin{figure}[ht]
\centering
\includegraphics[width=9cm]{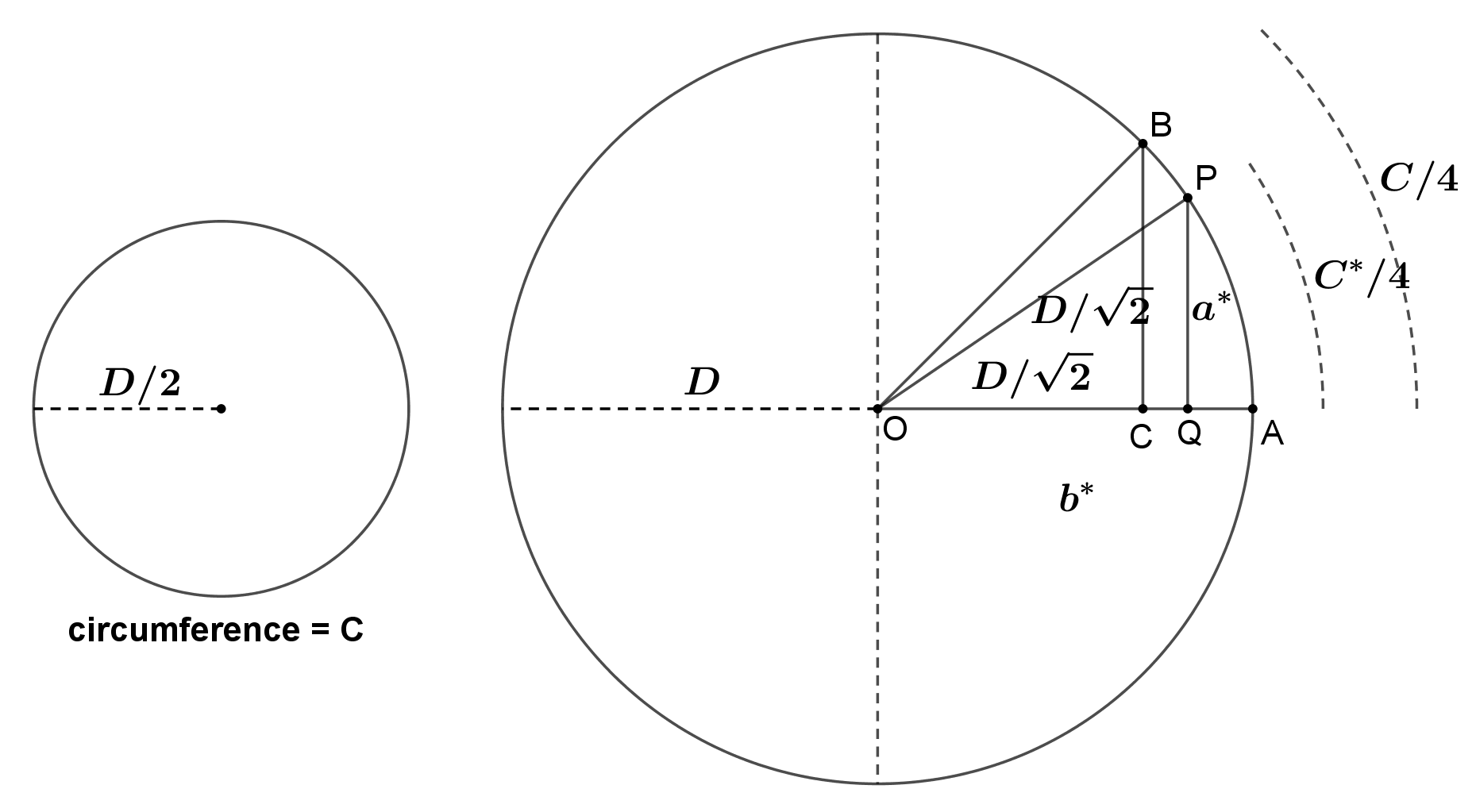}
\caption{Diagram for the description of the third method for the computation of the circumference of a circle}\label{ArcSinSeriesDerivation}
\end{figure}

Let $C$ be the circumference of a circle of diameter $D$. Let $C$ be unknown but an approximation $C^\ast$  to $C$ be known.   The procedure described in {\em Yuktibh\=a\d s\=a} yields a better approximation to $C$.

We consider a circle of radius $D$ and an arc of length $\frac{C^\ast}{4}$ on this circle (see Figure \ref{ArcSinSeriesDerivation}).
We calculate $a^\ast = \text{jy\=a}\,(\frac{C^\ast}{4})$ using the M\=adhava-Newton series for jy\=a:
\begin{equation}\label{aast}
a^\ast 
= \Big(\frac{C^\ast}{4}\Big) - \Big(\frac{1}{D}\Big)^2
\frac{1}{ 1 \cdot 2 \cdot 3}\Big(\frac{C^\ast}{4}\Big)^3 +
\Big(\frac{1}{D}\Big)^4\frac{1}{1 \cdot 2 \cdot 3 \cdot 4\cdot 5}\Big(\frac{C^\ast}{4}\Big)^5  -  \cdots. 
\end{equation}
We next calculate $b^\ast=\text{koj\=a}\,(\frac{a^\ast}{4})$ as follows:
$$
b^\ast = \sqrt{D^2-(a^\ast)^2}.
$$
The arc $\frac{C}{4}$ is one-eighth of the circle of radius $D$. Hence we have:\footnote{{\em Yuktibh\=a\d s\=a} has given the rationale of these results in an earlier section (see {\em Ga\d nita-Yukti-Bh\=a\d s\=a} Section 7.2.2; \cite{Sarma2008} p. 211).}
$$
\text{jy\=a}\,\Big(\frac{C}{4}\Big) = \sqrt{\frac{D^2}{2}}, \quad \text{kojy\=a}\,\Big(\frac{C}{4}\Big) = \sqrt{\frac{D^2}{2}}.
$$
We assume that $b^\ast > a^\ast$ so that $C>C^\ast$. By {\em j\=\i ve-parspara ny\=aya}, we have
\begin{align*}
\Delta & = \text{jy\=a}\,\Big(\frac{C}{4}-\frac{C^\ast}{4}\Big)\\
&
=
\frac{1}{D}\left[ \text{jy\=a}\,\Big(\frac{C}{4}\Big) 
\text{kojy\=a}\,\Big(\frac{C^\ast}{4}\Big) -
\text{kojy\=a}\,\Big(\frac{C}{4}\Big) \text{jy\=a}\,\Big(\frac{C^\ast}{4}\Big)\right] \\
&=
\frac{1}{D}\left[ \sqrt{\frac{D^2}{2}} b^\ast - \sqrt{\frac{D^2}{2}}a^\ast\right]\\
&=
\sqrt{\frac{(b^\ast)^2}{2}} -\sqrt{\frac{(a^\ast)^2}{2}}.
\end{align*}
Now, by Eq.\eqref{arcsin02}, we have
$$
\frac{C}{4}-\frac{C^\ast}{4}
\approx
\Delta + \frac{1}{D^2}\frac{\Delta^3}{6}\\
= \delta \quad\text{(say)}.
$$
Then
$$
C\approx C^\ast + 4\delta.
$$
This gives a better approximation to $C$ than $C^\ast$. If $b^\ast< a^\ast$, then $C<C^\ast$ and we will have to consider the difference $\frac{C^\ast}{4}-\frac{C}{4}$ in the above described procedure and $C$ will be given by $C=C^\ast -4\delta$. 
\renewcommand{\arraystretch}{1.35}
\begin{figure}[t]
\centering
\begin{tabular}{lcl}
$D$ & = & $1400^\prime$\\
$C^\ast$ & = & $D\times\frac{22}{7}=4400^\prime$\\
$\frac{C^\ast}{4}$ & = & $1100^\prime$\\
First term = $+ \tfrac{C^\ast}{4}$ & = & $+(1100^\prime)$\\
Second term = $ - (1100^\prime) (\tfrac{1}{D})^2
\tfrac{1}{2\cdot 3}(\frac{C^\ast}{4})^2 $ & = & $-(113^\prime\,\, 10^{\prime\prime}\,\,49^{\prime\prime\prime})$\\
Third term = $+(113^\prime\,\, 10^{\prime\prime}\,\,49^{\prime\prime\prime})(\tfrac{1}{D})^2\tfrac{1}{4\cdot 5}(\frac{C^\ast}{4})^2$& = & $+(3^\prime\,\, 29^{\prime\prime}\,\,37^{\prime\prime\prime})$\\
Fourth term = $-(3^\prime\,\, 29^{\prime\prime}\,\,37^{\prime\prime\prime})(\tfrac{1}{D})^2\tfrac{1}{6\cdot 7}(\tfrac{C^\ast}{4})^2$ & = & $-(0^\prime\,\, 3^{\prime\prime}\,\,5^{\prime\prime\prime})$\\
Fifth term = $+(0^\prime\,\, 3^{\prime\prime}\,\,5^{\prime\prime\prime})(\tfrac{1}{D})^2\tfrac{1}{8\cdot 9}(\tfrac{C^\ast}{4})^2$ & = & $+(0^\prime\,\, 0^{\prime\prime}\,\,2^{\prime\prime\prime})$\\
$a^\ast$ using Eq.\eqref{aast}  & = & 
$990^\prime\,\, 15^{\prime\prime}\,\, 45^{\prime\prime\prime}$\\
$(a^\ast)^2$ & = & $980619^\prime\,\, 49^{\prime\prime}\,\, 8^{\prime\prime\prime}$\\
$(b^\ast)^2 = D^2-(a^\ast)^2$ & = & $979380^\prime\,\, 10^{\prime\prime}\,\,52^{\prime\prime\prime}$\\
$\sqrt{\frac{(a^\ast)^2}{2}}$& = & $700^\prime\,\,13^{\prime\prime}\,\,17^{\prime\prime\prime}$\\
$\sqrt{\frac{(b^\ast)^2}{2}}$ & = & $699^\prime\,\,46^{\prime\prime}\,\,43^{\prime\prime\prime}$\\
Since $a^\ast>b^\ast$, we consider $\frac{C^\ast}{4}-\frac{C}{4}$. & & \\
$\Delta = \text{jy\=a}(\tfrac{C^\ast}{4}-\tfrac{C}{4}) = \sqrt{\frac{(a^\ast)}{2}} - \sqrt{\frac{(b^\ast)^2}{2}}$ & = & $0^\prime\,\, 26^{\prime\prime}\,\, 24^{\prime\prime\prime}$\\
$\delta = \Delta + \frac{1}{D^2}\frac{\Delta^3}{6}$ & = & $0^\prime\,\, 26^{\prime\prime}\,\, 24^{\prime\prime\prime}$\\
$4\delta$ & = & $1^\prime\,\, 45^{\prime\prime}\,\,36^{\prime\prime\prime}$\\
$C=D-4\delta$ & = & $4398^\prime\,\, 14^{\prime\prime}\,\,24^{\prime\prime\prime}$
\end{tabular} 
\caption{Computations} \label{Computations}
\end{figure}
\renewcommand{\arraystretch}{1}
\subsection{Illustrative example}
R\=amavarma (Maru) Tampur\=a has illustrated the procedure with the a numerical example. He has used the procedure to get the circumference $C$ of a circle of diameter $D=1400^\prime$. As an approximation to $C$, he has used the value of the circumference computed using the value $\frac{22}{7}$ for $\pi$. For details of the computations, see Figure \ref{Computations}.

The value of the circumference $C$ obtained here is still only an approximate value. The correct value of the circumference is %
$
3.14159265\times 1400^\prime = 4398^\prime\,\, 13^{\prime\prime}\,\,47^{\prime\prime\prime}.
$
%
%%%%%%%%%%%%%%%%%%%%%%%%%%%%%%%%%%%%%%%%%%%%%%%%%
%

%
\end{document}